\documentclass[aap,MSNbibl,citesort,seceqn,dvips]{arximspdf}

\doi{10.1214/10-AAP743}
\volume{21}
\issue{5}
\pubyear{2011}
\firstpage{1911}
\lastpage{1932}

\makeatletter

\newcommand{\rr}{\mathbb R}
\newcommand{\ev}{\mathbb E}
\newcommand{\nn}{\mathbb N}
\newcommand{\zz}{\mathbb Z}
\newcommand{\pp}{\mathbb P}

\newtheorem{prop}[theorem]{Proposition}
\newtheorem{theorem}{Theorem}[section]
\newtheorem{lemma}[theorem]{Lemma}
\newproclaim{assumption}[theorem]{Assumption}
\newproclaim{Remark}{Remark}

\makeatother

\begin{document}
\begin{frontmatter}

\title{Competing particle systems evolving by interacting L\'{e}vy
processes\thanksref{T1}}
\runtitle{Competing particle systems}

\thankstext{T1}{Supported in part by NSF Grant DMS-08-06211.}

\begin{aug}
\author[A]{\fnms{Mykhaylo} \snm{Shkolnikov}\corref{}\ead[label=e1]{mshkolni@math.stanford.edu}}
\runauthor{M. Shkolnikov}
\affiliation{Stanford University}
\address[A]{Department of Mathematics\\
Stanford University\\
Stanford, California 94305\\
USA\\
\printead{e1}} 
\end{aug}

\received{\smonth{2} \syear{2010}}
\revised{\smonth{6} \syear{2010}}

%
\begin{abstract}
We consider finite and infinite systems of particles on the real line
and half-line evolving in continuous time. Hereby, the particles are
driven by i.i.d. L\'{e}vy processes endowed with rank-dependent drift
and diffusion coefficients. In the finite systems we show that the
processes of gaps in the respective particle configurations possess
unique invariant distributions and prove the convergence of the gap
processes to the latter in the total variation distance, assuming a
bound on the jumps of the L\'{e}vy processes. In the infinite case we
show that the gap process of the particle system on the half-line is
tight for appropriate initial conditions and same drift and diffusion
coefficients for all particles. Applications of such processes include
the modeling of capital distributions among the ranked participants in
a financial market, the stability of certain stochastic queueing and
storage networks and the study of the Sherrington--Kirkpatrick model of
spin glasses.
\end{abstract}

%

\pdfkeywords{60J25, 60H10, 91B26, Stochastic differential equations,
Levy processes, semimartingale reflected Brownian motions,
Harris recurrence, capital distributions,
Levy queueing networks}

\end{frontmatter}

\section{Introduction}

Recently, invariant distributions for the gaps in a particle system on
the real line have received much attention. In the continuous time
setting such questions are motivated by the study of gaps in ordered
Brownian particle systems with rank-dependent drifts and diffusion
coefficients. The latter arise in the modeling of the capital
distribution in a financial market (see, e.g., \cite{bf,cp,fe}
and \cite{ip}) and as heavy traffic
approximations of queueing networks (see Section 5 of \cite{ip}
for the correspondence between the processes of gaps and reflected
Brownian motions and, e.g., \cite{ha,hn} and \cite{wi2} for
the heavy traffic approximation of queueing networks by the
latter). In the discrete time setting these questions appeared in the
study of the Sherrington--Kirkpatrick model of spin glasses and were
allowed to characterize the quasi-stationary states in the free energy
model (see \cite{ra} for this result and \cite{aa,sh} for its
extensions). In these papers the increments added to the
particles in each time step were either assumed to be i.i.d. or
normally distributed. The latter evolutions by Gaussian increments can
be viewed as a time discretized version of the Brownian particle
systems.

The present paper studies the processes of gaps in finite and infinite
particle systems on the real line and half-line evolving in continuous
time and in which particles are driven by i.i.d. L\'{e}vy processes
with jumps. In the finite case the latter are endowed with
rank-dependent drift and diffusion coefficients. We will refer to these
evolutions as competing particle systems. More precisely, the competing
particles are indexed by $i\in I$ with $I=\{1,\ldots,N\}$ in the case of
finitely many particles and $I=\nn$ in the case of infinitely many
particles. The evolution on the whole line will be referred to as the
unregulated evolution and the evolution on the half-line $[b,\infty)$
for some barrier $b\in\rr$ as the regulated evolution. The particle
systems are defined as the unique weak solutions (see Section \ref{sec21} and
Proposition~\ref{prop3.1} for more details) to the following stochastic
differential equations:
%
%
\begin{eqnarray}\label{sde1}
dX_i(t)&=&\sum_{j\in I} 1_{\{X_i(t)=X_{(j)}(t)\}}\delta_j \,dt +
\sum
_{j\in I} 1_{\{X_i(t)=X_{(j)}(t)\}}\sigma_j
\,dB_i(t)\nonumber\\[-8pt]\\[-8pt]
&&{} + dL_i(t) + dR_i(t) \nonumber
\end{eqnarray}
for all $i\in I$ where $R_i(t)\equiv0$, $i\in I$ in the unregulated
case and
%
%
\begin{equation}\label{sde2}\quad
R_i(t)=\sum_{0\leq s\leq t} 1_{\{X_i(s-)+\Delta
L_i(s)<b\}}\bigl(b-X_i(s-)-\Delta L_i(s)\bigr)+\Lambda_{(i,b)}(t),
\end{equation}
$i\in I$ in the regulated case. Hereby, $X_{(1)}(t)\leq X_{(2)}(t)\leq
\cdots$ is the ordered vector of particles in which ties are broken
according to an arbitrarily specified ordering of the initial
configuration, $B(t)=(B_i(t), i\in I)$ is a collection of i.i.d.
standard Brownian motions, $L(t)=(L_i(t), i\in I)$ is an independent
collection of i.i.d. pure jump L\'{e}vy processes each making finitely
many jumps on any finite time interval. For each $i\in I$ and $s\geq0$
the term $\Delta L_i(s)$ denotes the jump $L_i(s)-L_i(s-)$ of the
process $L_i$ at time $s$ and each $\Lambda_{(i,b)}(t)$ is the local
time process of $X_i$ at $b$. To avoid any confusion we choose the same
normalization of the local time processes as in \cite{bf,bg}
and \cite{ip} so that $\Lambda_{(i,b)}(t)$ is one
of the summands in the semi-martingale decomposition of the process
$\frac{1}{2}|X_i(t)-b|$ rather than $|X_i(t)-b|$. Moreover, we impose
the initial conditions $X_1(0)\leq X_2(0)\leq\cdots$ and
$0=R_1(0)=R_2(0)\cdots$ and, in addition, $X_1(0)\geq b$ in the
regulated case. To distinguish between the two evolutions we denote the
particle configuration at a time $t\geq0$ by $X(t)=(X_i(t)\dvtx i\in I)$
in the unregulated particle system and by $X^R(t)=(X^R_i(t)\dvtx i\in I)$
in the regulated particle system.

Heuristically, particle $i$ in the unregulated particle system is
driven by a L\'{e}vy process whose drift and diffusion coefficients
change according to its rank in the particle system. The regulated
particle system is defined similarly, except that whenever a particle
hits the barrier $b$ or jumps over $b$, it is reflected in the former
case and it is placed at $b$ in the latter case.

To simplify the notation we set $Y_i(t)=X_{(i)}(t)$ and
$Y^R_i(t)=X^R_{(i)}(t)$, respectively, for all $i\in I$. Next, we can
define the processes of gaps, which will be the main quantities of
interest, by
%
%
\begin{equation}
Z(t)\equiv(Z_1(t),Z_2(t),\ldots
)=\bigl(Y_2(t)-Y_1(t),Y_3(t)-Y_1(t),\ldots\bigr)
\end{equation}
and
%
%
\begin{equation}
Z^R(t)\equiv(Z^R_1(t),Z^R_2(t),\ldots)=\bigl(Y^R_1(t)-b,Y^R_2(t)-b,\ldots\bigr),
\end{equation}
respectively. For $I=\{1,\ldots,N\}$ we can set $C_i(t)=e^{X_i(t)}$ and
note that the properties of the gap process $Z(t)$ correspond to the
properties of the mass partition
%
%
\begin{equation}
\frac{C_{(i)}(t)}{C_{(1)}(t)+\cdots+C_{(N)}(t)}=\frac{1}{\sum_{j=1}^N
e^{X_{(j)}(t)-X_{(i)}(t)}},
\end{equation}
$1\leq i\leq N$ and the analogous statement is true for the process
$Z^R(t)$.

In models for the capital distribution in a financial market the
processes $C_1(t),\ldots,C_N(t)$ stand for the capitalizations of the
different firms and the particles represent the logarithmic
capitalizations of the firms. Due to this interpretation, we will call
the diffusion coefficients from now on volatilities or volatility
coefficients of the particles. In this context it is natural to allow
the particles to have jumps corresponding to such events as mergers of
firms, shares emissions or jumps in the stock price. Indeed, since
Merton's seminal paper \cite{me} and a vast number of works that
followed, it is widely agreed that market models allowing for
discontinuities in the stock price are able to better account
for sudden large changes in the values of the stocks, the
incompleteness of financial markets and the imperfectness of the
hedging strategies for options when compared to models assuming
continuous trajectories of the stock prices. Hereby, models in which
the jumps originate from pure jump L\'{e}vy processes play a central
role. For a recent summary of the research in this direction and a list
of references we refer to the book \cite{ct}.

The models on the half-line describe markets in which either firms are
bailed out if their capitalization falls below a certain predetermined
barrier or in which capitalizations are constrained not to overcome a
certain barrier, depending on the interpretation of the model. In the
absence of jumps such models can be analyzed by exploiting their
relation with semi-martingale reflected Brownian motions in the sense
of \cite{rw,tw} as was done in \cite{ip} for the
models on the whole line [see the proof of Lemma \ref{lemma2.1}(b) below for the
connection]. In \cite{ms} the authors considered a related model
for a market with two firms defaulting if their capitalization crosses
a fixed threshold. They were able to determine how to split the total
drift, representing the total amount of tax cuts or subsidies, between
the two firms to maximize the probability that both of them survive.

The mass partitions above correspond in this context to the market
weights of the ranked market participants. For this reason, all results
on the gap processes in the finite particle systems translate directly
into corresponding results on the processes of market weights. In
particular, the invariant distributions, the existence and uniqueness
of which we prove for the processes of gaps in the finite systems both
on the line and on the half-line, stand for capital distributions among
the ranked market participants which are not changing under the
evolution. We also show that the gap processes converge to the
respective unique invariant distributions in the total variation
distance. This corresponds to the statement that the considered
financial market approaches a stable capital allocation in a strong
sense. In addition, under the invariant distributions for the gap
processes, all gaps are finite by definition which implies that the
limiting stable capital distributions are nondegenerate, in the sense
that the market weights of all firms are positive almost surely.

In the context of stochastic queueing and storage networks, our model
is closely related to the so-called L\'{e}vy networks which are
described mathematically by L\'{e}vy processes regulated to stay in an
orthant by normal projections onto the boundary and normal reflections
at the boundary (see Chapter IX of \cite{as} and the references
therein for more details). In particular, the techniques we use to
prove that the processes of gaps possess unique invariant distributions
and converge to the latter can be used without alteration to show the
same statements for an integrable L\'{e}vy process with a nondegenerate
Brownian part, whose components are endowed with negative drifts and
which is regulated to stay in an orthant, provided that its jumps are
dominated by the drifts in an appropriate sense. The convergence in
total variation of the process to its invariant distribution proves the
stability of the corresponding queueing network, meaning that the joint
law of the workload processes converges in the strong sense.

Throughout the paper we make the following assumptions.
\begin{assumption}\label{assump1.1}
(a) In all cases we assume that the volatilities $\sigma_1,\sigma
_2,\ldots$ are positive and it holds that
%
%
\begin{equation}
\ev[|L_i(t)|]<\infty,\qquad \ev[L_i(t)]=0,\qquad i\in I.
\end{equation}

(b) If $I=\nn$, the drifts and volatilities are assumed to satisfy
%
%
\begin{equation}
\delta_M=\delta_{M+1}=\cdots,\qquad \sigma_1=\sigma_2=\cdots
\end{equation}
for some $M\geq1$ in the unregulated case and, in addition, $0>\delta
_1=\delta_2=\cdots$ in the regulated case.
\end{assumption}

Our results for the systems with finitely many particles are
summarized in the following two theorems.
\begin{theorem}[(Invariant distributions, finite systems)]\label{theo1.2}
\textup{(a)} If $I=\{1,\ldots,N\}$ and $L_1$ satisfies
%
%
\begin{equation}\label{cond1}
\ev\biggl[\sum_{0\leq s\leq1} |\Delta L_1(s)|\biggr]<\frac{1}{N}\cdot
\min
(\delta_1-\delta_2,\ldots,\delta_{N-1}-\delta_N),
\end{equation}
then the gap process $Z$ of the unregulated particle system has a
unique invariant distribution.

\textup{(b)} If $I=\{1,\ldots,N\}$ and $L_1$ satisfies
%
%
\begin{equation}\label{cond2}
\ev\biggl[\sum_{0\leq s\leq1} |\Delta L_1(s)|
\biggr]<\frac
{1}{N}\cdot\min(-\delta_1,\delta_1-\delta_2,\ldots,\delta
_{N-1}-\delta_N),
\end{equation}
then the gap process $Z^R$ of the regulated particle system has a
unique invariant distribution.
\end{theorem}
\begin{theorem}[(Convergence, finite systems)]\label{theo1.3}
\textup{(a)} If $I=\{1,\ldots,N\}$ and condition (\ref{cond1}) holds, then
the process $Z(t)$, $t\geq0$ converges in total variation to its unique
invariant distribution.

\textup{(b)} If $I=\{1,\ldots,N\}$ and condition (\ref{cond2}) holds, then
the same statement is true for the process $Z^R(t)$, $t\geq0$.
\end{theorem}

We remark that in the absence of jumps, condition (\ref
{cond1}) simplifies to $\delta_1>\cdots>\delta_N$. In particular, the
latter is stronger than condition (3.2) of \cite{ip} (used there
to establish the existence and uniqueness of an invariant distribution
for the gap process of the finite unregulated system without
jumps)\vspace*{1pt}
which in our notation reads $\frac{1}{i}\sum_{j=1}^i \delta_j-\frac
{1}{N}\sum_{j=1}^N \delta_j>0$, $1\leq i\leq N-1$. This shows that one
cannot expect the conditions in Theorems \ref{theo1.2} and \ref{theo1.3} to be sharp
in general.

For the infinite regulated system we show the following proposition.
\begin{prop}[(Tightness of the infinite regulated system)]\label{prop1.4}
If $I=\nn$ and
%
%
\begin{equation}\label{cond3}
\ev\biggl[\sum_{0\leq s\leq1} |\Delta L_1(s)|\biggr]<-\delta_1,
\end{equation}
then the family $Z^R(t)$, $t\geq0$ is tight on $\rr_+^\nn$ equipped
with the product topology for any initial condition in
%
%
\begin{equation}
W=\biggl\{0\leq z_1\leq z_2\leq\cdots\leq\big|\liminf_{i\rightarrow\infty}\frac
{z_i}{i}>0\biggr\}.
\end{equation}
\end{prop}

The proof of Proposition \ref{prop1.4} relies heavily on the fact that
under Assumption~\ref{assump1.1}, the evolution of each particle in the infinite
regulated system is given by a L\'{e}vy process with negative drift
regulated to stay on $[b,\infty)$. For this reason it does not carry
over to the setting of the infinite unregulated system. At this point
we also note that in the special case of the infinite regulated system
with equal drifts $\delta_1=\delta_2=\cdots,$ equal volatilities
$\sigma
_1=\sigma_2=\cdots$ and no jumps, the processes $X^R_i(t)-b$ are
independent reflected Brownian motions on $\rr_+$. Due to the
negativity of the drifts this implies that for each $i\in I$ the
process $X^R_i(t)-b$ converges in law to an exponential random
variable. Moreover, a sequence of i.i.d. exponential random variables
is almost surely dense in $\rr_+$ by the second Borel--Cantelli lemma.
This shows that, in general, one cannot expect the infinite regulated
system to have an invariant distribution, since already in the
described special case the number of particles on each nonempty
interval of the form $[b,y)$ will tend to infinity for any initial
particle configuration. In contrast to this, in \cite{ra} the
authors show the existence of an infinite family of invariant
distributions for the gap process in the corresponding particle system
on the whole line. In the more general case of nonconstant drifts or
volatilities the questions of existence and uniqueness of invariant
distributions for the gap process of the infinite regulated system and
of the convergence of the latter to an invariant distribution are open.
The approach used in the analysis of the finite systems may apply to
this case as well. However, one has to establish the recurrence and
irreducibility structure of the gap process in an infinite-dimensional
space where the analysis of the corresponding reflected Brownian motion
is more intricate.

In \cite{pp} the authors treat the infinite unregulated system
without jumps with the drift sequence $\delta_1,0,0,\ldots$ and a
constant and positive volatility sequence $\sigma_1=\sigma_2=\cdots$ and
are able to find an invariant distribution for the gap process on a
subset of $W$. In \cite{ra} the authors characterize all invariant
distributions for the process of gaps in the case of constant drift and
volatility sequences. As in the regulated case the questions of
existence and uniqueness of invariant distributions and of the
convergence of the gap process to the latter for a general drift or
volatility sequence are open. The main obstacle hereby is the same as
in the case of the infinite regulated system (see the end of the
previous paragraph).

It is a natural question to ask to provide an explicit description of
the invariant distributions in Theorem \ref{theo1.2} in the presence of jumps.
However, even in the case of the regulated particle system with a
single particle, where the gap process is a one-dimensional L\'{e}vy
process regulated to stay nonnegative, the invariant distribution is
not known explicitly. To the best of our knowledge the most explicit
result in this direction is obtained in \cite{ba}, Section 3.
There, in the described special case with the additional assumption
that the magnitude of the negative jumps of the driving L\'{e}vy
process is exponentially distributed, an expression for the Fourier
transform of the invariant distribution is given in terms of the
characteristic exponent of the driving L\'{e}vy process and the
expected accumulated local time of the gap process at zero over a unit
time interval when the initial value of the gap process is chosen
according to the unique invariant distribution (see \cite{ba},
equation (3.11)). A completely explicit description of the invariant
distribution for the gap process, even for specific examples of a
regulated or an unregulated system with multiple particles in which
jumps are present, seems to be out of reach at the moment.

Previous papers on Brownian systems with rank-dependent drifts and
volatilities use extensively the results of \cite{wi} on reflected
Brownian motions and classical results on constrained diffusion
processes to obtain the existence and uniqueness of the invariant
distribution and the convergence of the gap process to the latter. Due
to the presence of jumps, these tools do not apply here and are
replaced by more general techniques from the ergodic theory of Harris
recurrent Markov processes, as in \cite{mt}. More precisely, to
obtain the results of Theorems \ref{theo1.2} and~\ref{theo1.3}, we first find the
dynamics of the gap processes in Lemma \ref{lemma2.1} below. We use it
subsequently to prove the tightness of the processes of gaps by
dominating each of them by an appropriate Harris recurrent Markov
process. Next, we show the Feller property for the gap processes to
establish the existence of the invariant distributions for the latter.
Finally, to prove the uniqueness of the invariant distributions and the
convergence of the processes of gaps, we employ the general theory of
Harris recurrent Markov processes.

The paper is structured as follows. In Section \ref{sec2} we treat the particle
systems with finitely many particles. We show first their existence and
explain some of their properties in Section \ref{sec21} and then prove
Theorems \ref{theo1.2} and \ref{theo1.3} in Section~\ref{sec22}.
Section \ref{sec3} deals with systems of
infinitely many particles. We show that the latter exist and are well
defined for appropriate initial conditions in Section \ref{sec31} and prove
Proposition \ref{prop1.4} in Section \ref{sec32}.

\section{Systems of finitely many particles}\label{sec2}

\subsection{Existence and properties of the processes}\label{sec21}

Throughout this section we deal with the two evolutions of finitely
many particles, that is, we set $I=\{1,\ldots,N\}$. The existence and
uniqueness of a weak solution to the unregulated version of (\ref
{sde1}) can then be seen as follows. As remarked in Section 2 of
\cite{bf}, the main result of \cite{bp} implies that for any
$x\in\rr^N$ with $x_1\leq\cdots\leq x_N$ there exists a unique weak
solution to
%
%
\begin{eqnarray}\quad
\label{sdec1}
dX^{c,x}_i(t)&=&\sum_{j=1}^N 1_{\{
X^{c,x}_i(t)=X^{c,x}_{(j)}(t)\}}\delta_j \,dt
+ \sum_{j=1}^N 1_{\{X^{c,x}_i(t)=X^{c,x}_{(j)}(t)\}}\sigma_j
\,dB_i(t),\\
\label{sdec2}
X^{c,x}(0)&=&x
\end{eqnarray}
defined on some probability space $(\Omega^x,\mathcal{F}^x,\pp^x)$.
Next, let $(\Omega^{x,l},\mathcal{F}^{x,l},\pp^{x,l})$, $l\in\nn$, be
copies of $(\Omega^x,\mathcal{F}^x,\pp^x)$ such that on each of them a
process $X^{c,x,l}$ of the same law as $X^{c,x}$ is defined. Moreover,
let $(\Omega^L,\mathcal{F}^L,\pp^L)$ be the probability space on which
the processes $L_1(t),\ldots,L_N(t)$ are defined. Then a weak solution
to (\ref{sde1}) can be defined on the product space
%
%
\begin{equation}
(\Omega,\mathcal{F},\pp)
=\biggl(\prod_{x,l}\Omega^{x,l}\times\Omega^L,\bigotimes
_{x,l}\mathcal
{F}^{x,l}\otimes\mathcal{F}^L,\bigotimes_{x,l} \pp^{x,l}\otimes\pp
^L\biggr).
\end{equation}
Denoting by $T_1<T_2<\cdots$ the jump times of the $\rr^N$-valued
process $L(t)=(L_1(t),\ldots,L_N(t))$ we set
%
%
\begin{eqnarray}
X(0)&=&(X_1(0),\ldots,X_N(0)),\\
X(t)&=&X^{c,X(T_k),k}(t-T_k),\qquad T_k<t<T_{k+1}, k\geq0,\\
X(T_k)&=&X(T_k-)+\Delta L(T_k),\qquad k\geq1,
\end{eqnarray}
where we have used the notation $T_0=0$. This gives a weak solution to
the unregulated version of (\ref{sde1}) by construction. Noting that
for any other weak solution the joint distribution of its continuous
part, its jump times and its jump sizes have to coincide with the
corresponding quantity of the solution just constructed, we conclude
that the weak solution is unique. Here we have used the uniqueness of
the weak solution to (\ref{sdec1}), (\ref{sdec2}).

The regulated process $X^R$ can be constructed on the same probability
space as $X$ as follows. We start with the desired initial value
$X^R(0)$ and set
\begin{eqnarray*}
X^R_i(t)&=&X^{c,X^R(T_k),k}_i(t-T_k)+\Lambda^{c,T_k}_{(i,b)}(t),\qquad
T_k<t<T_{k+1}, k\geq0,\\
X^R_i(T_k)&=&
\cases{X^R_i(T_k-)+\Delta L_i(T_k), &\quad if
$X^R_i(T_k-)+\Delta L_i(T_k)\geq b$,\cr
b, &\quad if $X^R_i(T_k-)+\Delta L_i(T_k)<b$,}\qquad k\geq1,
\end{eqnarray*}
for all\vspace*{1pt} $1\leq i\leq N$ where $\Lambda^{c,T_k}_{(i,b)}(t)$, $t\geq
T_k$,
is the local time at $b$ of the process $X^R_i(t)$, $t\geq T_k$. We
remark that this is precisely the construction of $X$ with the
regulation and reflection at the barrier being added at the appropriate
random times. Due to the definition of the processes $X^{c,x,l}$ and of
the involved local time processes, $X^R$ and the corresponding process
of regulations solve the regulated system (\ref{sde1}), (\ref{sde2}).
The uniqueness of the weak solution to (\ref{sdec1}), (\ref{sdec2})
implies that the weak solution of the regulated system is unique.

It follows immediately that the gap processes $Z(t)$, $Z^R(t)$ are
well defined and unique in law. Their dynamics are given in the next lemma.
\begin{lemma}\label{lemma2.1}
\textup{(a)} The components of the gap process $Z$ in the unregulated particle
system satisfy
\begin{eqnarray*}
dZ_i(t)&=&(\delta_{i+1}-\delta_1)\,dt + \sigma_{i+1} \,d\beta
_{i+1}(t)-\sigma
_1 \,d\beta_1(t)+d\lambda_{i+1}(t)-d\lambda_1(t)\\
&&{}+\tfrac{1}{2}\bigl(d\Lambda_{(i,i+1)}(t)-d\Lambda
_{(i+1,i+2)}(t)+d\Lambda
_{(1,2)}(t)\bigr)\\
&&{}+\bigl(F_i(Z(t-),\Delta\lambda(t))-\bigl(Z_i(t-)+\Delta\bigl(\lambda
_{i+1}(t)-\lambda_1(t)\bigr)\bigr)\bigr).
\end{eqnarray*}

\textup{(b)} The components of the gap process $Z^R$ in the regulated particle
system are governed by
\begin{eqnarray*}
dZ^R_i(t)&=&\delta_i \,dt + \sigma_i \,d\beta_i(t)+d\lambda_i(t)+\tfrac
{1}{2}\bigl(d\Lambda^R_{(i-1,i)}(t)-d\Lambda^R_{(i,i+1)}(t)\bigr)
1_{\{
i\neq1\}}\\
&&{}+\bigl(d\Lambda^R_{(0,1)}(t)-\tfrac{1}{2}\,d\Lambda^R_{(1,2)}(t)\bigr)
1_{\{
i=1\}}\\
&&{}+\bigl(F^R_i(Z^R(t-),\Delta\lambda(t))-\bigl(Z^R_i(t-)+\Delta\lambda
_i(t)\bigr)\bigr).
\end{eqnarray*}
Hereby, $\beta_1(t),\ldots,\beta_N(t)$ are i.i.d. standard
Brownian motions, $\lambda(t)=(\lambda_1(t),\ldots,\break\lambda_N(t))$
is a
L\'{e}vy process whose components are i.i.d. pure jump L\'{e}vy
processes with the same distribution as $L_1(t)$ and independent of
$\beta_1(t),\ldots,\beta_N(t)$,\break $\Lambda_{(i,i+1)}(t)$ and $\Lambda
^R_{(i,i+1)}(t)$ are the local times at $0$ of the processes
$Y_{i+1}(t)-Y_i(t)$ and $Y^R_{i+1}(t)-Y^R_i(t)$, respectively, with
$\Lambda_{(N,N+1)}(t)\equiv\Lambda^R_{(N,N+1)}(t)\equiv0$, and
$\Lambda
^R_{(0,1)}(t)$ is the local time of the process $Y^R_1(t)$ at $b$.
Finally, $F$ and $F^R$ describe the value of the gap process after a
jump as a function of its value before the jump and the jump of
$\lambda
$. More explicitly,
\begin{eqnarray*}
F_i(z,\eta)&=&\min_{1\leq j_1<\cdots<j_{i+1}\leq N} \max
(z_{j_1-1}+\eta
_{j_1},\ldots,z_{j_{i+1}-1}+\eta_{j_{i+1}})\\
&&{}-\min_{1\leq j\leq N} (z_{j-1}+\eta_j)
\end{eqnarray*}
for all $1\leq i\leq N-1$, $z=(z_1,\ldots,z_{N-1})\in\rr_+^{N-1}$ and
$\eta=(\eta_1,\ldots,\eta_N)\in\rr^N$ with the convention $z_0=0$ and
\[
F^R_i(z,\eta)=\min_{1\leq j_1<\cdots<j_i\leq N} \max\bigl(\max
(z_{j_1}+\eta_{j_1},0),\ldots,\max(z_{j_i}+\eta_{j_i},0)\bigr)
\]
for all $1\leq i\leq N$, $z=(z_1,\ldots,z_N)\in\rr_+^N$ and $\eta
=(\eta
_1,\ldots,\eta_N)\in\rr^N$. Moreover, the state spaces the processes are
given by
\[
W^{N-1}=\{z_1,\ldots,z_{N-1}|0\leq z_1\leq z_2\leq\cdots\leq z_{N-1}\}
\subset\rr_+^{N-1}
\]
and $W^N$, respectively.
\end{lemma}
\begin{Remark*} The contribution of the terms
\[
\bigl(F_i(Z(t-),\Delta\lambda(t))-\bigl(Z_i(t-)+\Delta\bigl(\lambda
_{i+1}(t)-\lambda_1(t)\bigr)\bigr)\bigr),\qquad 1\leq i\leq N-1,
\]
and
\[
\bigl(F^R_i(Z^R(t-),\Delta\lambda(t))-\bigl(Z^R_i(t-)+\Delta\lambda
_i(t)\bigr)
\bigr),\qquad 1\leq i\leq N,
\]
in the respective dynamics can be understood as follows. The
contribution is nonzero if and only if one of the particles jumps and
this jump changes the ranks of the particles. If this jump does not
change the left-most particle in the unregulated system or does not
involve a regulation in the regulated system, then these terms
correspond to consecutive normal reflections of the process of gaps at
faces of $W^{N-1}$ or $W^N$, respectively. More precisely, the gap
process is normally reflected at the faces $\{z_{j_1}=z_{j_1+1}\},\ldots
,\{z_{j_2-1}=z_{j_2}\}$, if the inequalities $z_{j_1}\leq
z_{j_1+1},\ldots,z_{j_1}\leq z_{j_2}$ are violated by the jump and at
the faces $\{z_{j_2}=z_{j_2-1}\},\ldots,\{z_{j_1+1}=z_{j_1}\}$, if the
inequalities $z_{j_1}\leq z_{j_2},\ldots,z_{j_2-1}\leq z_{j_2}$ are
violated by the jump. If the jump changes the left-most particle in the
unregulated system, then the particles are relabeled and this term
gives the change of the gaps due to relabeling. If a particle jumps
below the barrier in the regulated evolution, then the particle
configuration is regulated and the particles are relabeled. The term
then gives the change of gaps due to both these operations.
\end{Remark*}
\begin{pf*}{Proof of Lemma \ref{lemma2.1}}
(a) For any fixed $t\geq0$ let $\pi_t\dvtx\{
1,\ldots,N\}\rightarrow\{1,\ldots,N\}$ be a bijection such that
%
%
\begin{equation}
X_{\pi_t^{-1}(1)}(t-)\leq X_{\pi_t^{-1}(2)}(t-)\leq\cdots\leq X_{\pi
_t^{-1}(N)}(t-).
\end{equation}
In Section 3 of \cite{bf} the authors show that in the absence of
jumps, that is, when $L(t)\equiv0$, it holds
%
%
\begin{equation}\label{sdey}
dY_i(t)=\sum_{j=1}^N 1_{\{\pi_t(j)=i\}}\,dX_j(t)+\frac
{1}{2}\bigl(d\Lambda_{(i-1,i)}(t)-d\Lambda_{(i,i+1)}(t)\bigr),
\end{equation}
where $\Lambda_{(i,i+1)}(t)$ are defined as in the statement of the
lemma for $1\leq i\leq N$ and we have set $\Lambda_{(0,1)}(t)\equiv0$.
Moreover, equation (3.3) of \cite{bf} states that there exist
i.i.d. standard Brownian motions $\beta_1(t),\ldots,\beta_N(t)$ such that
%
%
\begin{equation}
dY_i(t)=\delta_i \,dt + \sigma_i \,d\beta_i(t)+\tfrac{1}{2}
\bigl(d\Lambda
_{(i-1,i)}(t)-d\Lambda_{(i,i+1)}(t)\bigr).
\end{equation}
This yields immediately the claim of part (a) of the lemma if
$L(t)\equiv0$. In the presence of jumps we define the pure jump
processes $\lambda_1(t),\ldots,\lambda_N(t)$ by $\Delta\lambda
_i(t)=\Delta L_{\pi_t^{-1}(i)}(t)$. The latter are i.i.d. pure jump L\'
{e}vy processes with the same law as $L_1(t)$. Indeed, the jump times
of $\lambda(t)$ and $L(t)$ coincide and the law of $\Delta\lambda(t)$
conditional on $\Delta\lambda(t)\neq0$ is the same as the law of
$\Delta L(t)$ conditional on $\Delta L(t)\neq0$. Next, we recall that
$T_k$, $k\geq0$ were defined by $T_0=0$ and as the jump times of the
process $L(t)$ for $k\geq1$. Since $X(t)$ coincides with the
corresponding process in the absence of jumps for $T_{k-1}<t<T_k$ and
any $k\geq1$ by its construction, it remains to verify that $\Delta
Z_i(T_k)$ coincides with the jump given by the right-hand side of the
equation in part (a) of the lemma for all $1\leq i\leq N-1$ and $k\geq
1$. But this follows directly from the definition of $F$.

(b) As in (a), we first treat the case $L(t)\equiv0$. In this case we
will show that the processes $\widetilde
{Z}_i^R(t)=Y^R_i(t)-Y^R_{i-1}(t)$, $1\leq i\leq N$, with $Y^R_0(t)\equiv
b$ follow dynamics corresponding to the dynamics in part (b) of the
lemma. To this end, we set
%
%
\begin{equation}
N_j(t)=\bigl|\bigl\{1\leq i\leq N| X^R_i(t)=X^R_{(j)}(t)\bigr\}\bigr|
\end{equation}
for all $1\leq j\leq N$, $t\geq0$. Next, we observe that Theorem 2.3 of
\cite{bg} applied to the continuous semi-martingales
$X^R_1(t),\ldots,X^R_N(t)$ yields for $1\leq j\leq N$,
\begin{eqnarray*}
dY^R_j(t)&=&\sum_{i=1}^N (N_j(t))^{-1} 1_{\{X^R_{(j)}(t)=X^R_i(t)\}}
\,dX^R_i(t)+(N_j(t))^{-1}\sum_{k=1}^{j-1} d\Lambda^R_{(k,j)}(t)\\
&&{}
-(N_j(t))^{-1}\sum_{k=j+1}^N d\Lambda^R_{(j,k)}(t),
\end{eqnarray*}
where $\Lambda^R_{(j_1,j_2)}(t)$\vspace*{1pt} is the local time of
$Y^R_{j_2}(t)-Y^R_{j_1}(t)$ at zero for $1\leq j_1<j_2\leq N$. Plugging
(\ref{sde1}) and (\ref{sde2}) into the latter equation and applying the
strong Markov property of $(\widetilde{Z}^R_1(t),\ldots,\widetilde
{Z}^R_N(t))$ to the entrance times of the set
\[
\partial_\varepsilon\rr_+^N=\{z\in\rr_+^N| \mathrm{dist}(z,\partial\rr
_+^N)\geq
\varepsilon\}
\]
for a fixed $\varepsilon>0$, one shows that $\widetilde
{Z}^R(t)=(\widetilde
{Z}^R_1(t),\ldots,\widetilde{Z}^R_N(t))$ evolves as a Brownian motion
\[
\bigl(\sigma_1\beta_1(t),\sigma_2\beta_2(t)-\sigma_1\beta_1(t),\ldots
,\sigma
_N\beta_N(t)-\sigma_{N-1}\beta_{N-1}(t)\bigr)
\]
with constant drift vector $(\delta_1,\delta_2-\delta_1,\ldots
,\delta
_N-\delta_{N-1})$ between an entrance time of the set $\partial
_\varepsilon
\rr_+^N$ and the first hitting time of $\partial\rr_+^N$ after that.
This is due to an application of Knight's theorem in the form of
\cite{ry}, page 183, to the martingale parts of $Y^R_j(t)$, $1\leq
j\leq N$. Consequently, letting $\varepsilon$ tend to zero, we observe
that $\widetilde{Z}^R(t)$ is a semi-martingale reflected Brownian
motion in the orthant $\rr_+^N$ in the sense of \cite{rw,tw}.
Now, Lemma 2.1 of \cite{tw} shows that the
Lebesgue measure of the set $\{t\geq0| \widetilde{Z}^R(t)\in
\partial\rr
_+^N\}$ is zero almost surely and Theorem 1 of \cite{rw} implies
that the times for which either $N_j(t)\geq3$ or $N_j(t)=2$ and
$\widetilde{Z}_1(t)=0$ do not contribute to the dynamics of $Y^R_j(t)$
for all $1\leq j\leq N$. Hence, the dynamics simplifies to
\begin{eqnarray*}
dY^R_j(t)&=&\delta_j \,dt + \sigma_j \,d\beta_j(t) + \bigl(d\Lambda
^R_{(0,1)}(t)-\tfrac{1}{2}\,d\Lambda^R_{(1,2)}(t)\bigr)\cdot1_{\{j=1\}
}\\
&&{}+\tfrac{1}{2}\bigl(d\Lambda^R_{(j-1,j)}(t)-d\Lambda
^R_{(j,j+1)}(t)
\bigr)\cdot1_{\{j\neq1\}}
\end{eqnarray*}
for all $1\leq j\leq N$. This yields immediately the statement of part
(b) of the lemma in the absence of jumps. The general case follows by
defining $\lambda(t)$ in the same way as in the proof of part (a) of
the lemma and by making the next two observations. First, the dynamics
of the process $Z^R(t)$ between the jump times of the process $L(t)$
coincides with the dynamics of the system in the absence of jumps as a
consequence of the construction of the process $X^R(t)$. Second, the
jumps of the process $Z^R(t)$ coincide with the jumps of the right-hand
side of the equation in part (b) of the lemma due to the definition of
$F^R$.
\end{pf*}

\subsection{Invariant distributions and convergence}\label{sec22}

In this section we investigate the existence and uniqueness of
invariant distributions of the gap processes in the two finite particle
systems, as well as the convergence of the processes of gaps to the
respective invariant distributions. We start with the proof of Theorem
\ref{theo1.2}.
\begin{pf*}{Proof of Theorem \ref{theo1.2}}
(a) (1) We first prove that the family
$Z(t)$, $t\geq0$, is tight for any initial value $z$ in $W^{N-1}$. To do
this, it suffices to show that the family $\widetilde{Z}_i(t)\equiv
Z_i(t)-Z_{i-1}(t)$, $t\geq0$, is tight for all $1\leq i\leq N-1$ where
we have set $Z_0(t)\equiv0$. To this end, we fix a $1\leq i\leq N-1$,
set $y=\widetilde{Z}_i(0)$ and define the process $D^y_i$ on the same
probability space as $X$ by
%
%
\begin{eqnarray}
D^y_i(t)&=&y+(\delta_{i+1}-\delta_i)t+\sigma_{i+1}\beta
_{i+1}(t)-\sigma
_i\beta_i(t)\nonumber\\[-8pt]\\[-8pt]
&&{}+\sum_{j=1}^N\sum_{0\leq s\leq t}|\Delta\lambda_j(s)|+\Lambda
^{D^y_i}(t), \nonumber
\end{eqnarray}
where $\Lambda^{D^y_i}(t)$ is the local time at $0$ of the process
$D^y_i(t)$ and $\beta_1(t),\ldots,\beta_N(t)$ and $\lambda
_1(t),\ldots
,\lambda_N(t)$ are the same as in the dynamics of $Z(t)$ given in
Lem\-ma~\ref{lemma2.1}(a). We note that $\Delta\widetilde{Z}_i(t)\leq\Delta D_i^y(t)$ for
all $t\geq0$ and that after each time $t\geq0$ with $\widetilde
{Z}_i(t)=D^y_i(t)$ the processes $\widetilde{Z}_i$ and $D^y_i$ evolve
in the same way until either the $(i-1)$st ranked particle and the
$i$th ranked particle collide, or the $(i+1)$st ranked particle and the
$(i+2)$nd ranked particle collide, or there is a jump of $(\lambda
_1(t),\ldots,\lambda_N(t))$. In particular, for any $t\geq0$ the
accumulated local time at $0$ on the set $\{0\leq s\leq t|\widetilde
{Z}_i(s)=D^y_i(s)=0\}$ is the same for the processes $\widetilde{Z}_i$
and~$D^y_i$. Putting these observations together, we conclude that
$\widetilde{Z}_i(t)\leq D^y_i(t)$ almost surely for all $t\geq0$.
Hence, to prove that the family $\widetilde{Z}_i(t)$, $t\geq0$ is
tight, it suffices to show that the family $D^y_i(t)$, $t\geq0$ is
tight.

(2) Next, we fix an $\varepsilon>0$ and find a $C=C(\varepsilon)>0$
such that
%
%
\begin{equation}\label{trap}
\pp\Bigl(\sup_{0\leq u\leq\varepsilon} D^0_i(u)\leq C\Bigr)>0.
\end{equation}
We claim that the process $D^y_i(n\varepsilon)$, $n\in\nn$, is a recurrent
Harris chain on $\rr_+$ with respect to the set $[0,C]$ in the sense of
Section 5.6 of \cite{du}. Indeed, the Harris property follows
immediately from the Harris property of the corresponding process in
the absence of jumps. Moreover, the law of large numbers for L\'{e}vy
processes and condition (\ref{cond1}) show that almost surely
\begin{eqnarray*}
&&\lim_{t\rightarrow\infty}
\frac{1}{t}\Biggl((\delta_{i+1}-\delta_i)t+\sigma_{i+1}\beta
_{i+1}(t)-\sigma_i\beta_i(t)+\sum_{j=1}^N\sum_{0\leq s\leq
t}|\Delta
\lambda_j(s)|\Biggr)\\
&&\qquad=\delta_{i+1}-\delta_i+\sum_{j=1}^N \ev\biggl[\sum_{0\leq s\leq
1}|\Delta\lambda_j(s)|\biggr]<0.
\end{eqnarray*}
Hence, there exist stopping times $0<\tau_1<\tau_2<\cdots$ tending to
infinity with $D^y_i(\tau_k)=0$ for all $k\in\nn$ with probability $1$.
This follows by setting
\begin{eqnarray*}
\tau_1&=&\inf\{u\geq1|D^y_i(u)=0\},\\
\tau_{k+1}&=&\inf\{u\geq\tau_k+1|D^y_i(u)=0\},\qquad k\geq1,
\end{eqnarray*}
and coupling $D_i^y(t)$, $t\geq s$ with the L\'{e}vy process
\begin{eqnarray*}
&&D_i^y(s)+(\delta_{i+1}-\delta_i)(t-s)+\sigma_{i+1}\bigl(\beta
_{i+1}(t)-\beta
_{i+1}(s)\bigr)-\sigma_i\bigl(\beta_i(t)-\beta_i(s)\bigr)\\
&&\qquad{}+\sum_{j=1}^N\sum_{s<u\leq t}|\Delta\lambda_j(u)|,\qquad t\geq
s,
\end{eqnarray*}
to conclude
\[
\pp\bigl(\{D_i^y(s)=0\}\cup\{\exists t>s|D_i^y(t)=0\}\bigr)=1
\]
for all $s\geq0$. In addition, we let $\tau_k(\varepsilon)$ be the integer
multiple of $\varepsilon$ which is closest to $\tau_k$ from above for all
$k\in\nn$ and observe that the strong Markov property of $D_i$ applied
to the stopping times $\tau_k$, $k\in\nn$, (\ref{trap}) and the second
Borel--Cantelli lemma imply that $D^y_i(\tau_k(\varepsilon))\leq C$ for
infinitely many $k\in\nn$ almost surely. This shows the recurrence of
$D^y_i(n\varepsilon)$, $n\in\nn$. Moreover, by the results of Section 5.6c
of \cite{du}, the chain is aperiodic and converges in total
variation to its unique invariant distribution which we denote by $\nu
_\varepsilon$. Furthermore, the uniqueness of invariant distributions of
recurrent Harris chains implies that $\nu_1=\nu_{{1}/{2}}=\nu
_{{1}/{4}}=\cdots.$ Next, we fix a $\zeta>0$ and a $w>0$ and define the
function $f_w\dvtx\rr_+\rightarrow\rr_+$ by $f_w(v)=e^{-w\cdot v}$. Since
$D_i$ is a Feller process [this can be shown along the lines of step
(3) below], its semi-group of transition operators is strongly
continuous on the space of continuous functions on $\rr_+$ vanishing at
infinity (see, e.g., Theorem 17.6 in \cite{ka}). Thus, denoting the
semi-group of transition operators corresponding to $D_i(t)$, $t\geq0$,
by $P^{D_i}(t)$, $t\geq0$, we can find a $Q\in\nn$ such that
%
%
\begin{equation}
\forall 0\leq t\leq2^{-Q}\dvtx \|P^{D_i}(t)f_w-f_w\|_\infty<\zeta.
\end{equation}
Moreover, for each $t\geq0$ we let $n(t)$ be the largest integer such
that $n(t)\cdot2^{-Q}\leq t$. All in all, setting $\mu_{y,t}=\delta_y
P^{D_i}(t)$ and writing $\mu(f)$ for $\int f \,d\mu$, we get for any
$t\geq0$,
\begin{eqnarray*}
&&|\mu_{y,t}(f_w)-\nu_1(f_w)|\\
&&\qquad\leq\bigl|\mu_{y,t}(f_w)-\mu_{y,n(t)
2^{-Q}}(f_w)\bigr|+\bigl|\mu_{y,n(t)2^{-Q}}(f_w)-\nu_1(f_w)\bigr|\\
&&\qquad=\bigl|\mu_{y,n(t)2^{-Q}}\bigl(P^{D_i}\bigl(t-n(t)2^{-Q}\bigr)f_w-f_w\bigr)\bigr|+\bigl|\mu
_{y,n(t)2^{-Q}}(f_w)-\nu_1(f_w)\bigr|\\
&&\qquad\leq\zeta+\bigl|\mu_{y,n(t)2^{-Q}}(f_w)-\nu_1(f_w)\bigr|.
\end{eqnarray*}
Taking first the limit $t\rightarrow\infty$ and then the limit $\zeta
\downarrow0$, we conclude that $\mu_{y,t}(f_w)$ converges to $\nu
_1(f_w)$. Since $w>0$ was arbitrary, it follows that the Laplace
transforms of the measures $\delta_y P^{D_i}(t)$ converge point-wise to
the Laplace transform of~$\nu_1$, hence, the measures $\delta_y
P^{D_i}(t)$ converge weakly to $\nu_1$. Thus, the family $D^y_i(t)$,
$t\geq0$, converges in law to $\nu_1$ and is, therefore, tight as
claimed.

(3) Next, denote by $P(t)$, $t\geq0$, the Markov semi-group of operators
corresponding to the Markov process $Z$. The tightness of $Z(t)$,
$t\geq
0$ implies that the family of probability measures $\frac{1}{t}\int_0^t
\delta_z P(s) \,ds$, $t>0$, is tight. Hence, we can find a sequence
$t_1<t_2<\cdots$ tending to infinity such that the weak limit
%
%
\begin{equation}
\nu\equiv\lim_{n\rightarrow\infty} \frac{1}{t_n}\int_0^{t_n}
(\delta_z
P(s)) \,ds
\end{equation}
exists and is a probability measure on $W^{N-1}$. We claim that $\nu$
is an invariant distribution of $Z$. To this end, we first remark that
$Z$ is a Feller process. For the evolution without jumps this is a
consequence of the Feller property of the reflected Brownian motion
(see Theorem 1.1 and the following remark in \cite{wi}) and the
observation that the process $(\widetilde{Z}_1(t),\ldots,\widetilde
{Z}_{N-1}(t))$ can be viewed as a reflected Brownian motion in the
sense of \cite{wi} (see \cite{ip}, Section 5). In our case
the Feller property can be seen as follows. Let $z^0$ be an arbitrary
point in $W^{N-1}$ and $(z^n)_{n=1}^\infty$ be a sequence in $W^{N-1}$
converging to it. Moreover, let $Z^{z^0}$ and $Z^{z^n}$, $n\geq1$, be
the gap processes with initial values $z^0$ and $z^n$, $n\geq1$,
respectively. We need to show that for any fixed $t\geq0$ the random
vectors $Z^{z^n}(t)$ converge in law to $Z^{z^0}(t)$. To this end, we
note that the gap process $Z$ can be constructed by first generating
the sequence of its jump times $T_1<T_2<\cdots$ and the corresponding
jumps $\Delta\lambda(T_1),\Delta\lambda(T_2),\ldots$ of $\lambda
(t)=(\lambda_1(t),\ldots,\lambda_N(t))$ and then defining $Z$
conditional on these choices by
%
%
\begin{equation}\label{repr}
Z(t)=\sum_{k=0}^\infty1_{\{T_k\leq t<T_{k+1}\}} Z^{c,F(Z(T_k-),\Delta
\lambda(T_k)),k}(t-T_k),
\end{equation}
where we have set $T_0=0$, $Z(0-)=Z(0)$, $F$ is the continuous function
as in part (a) of the Lemma \ref{lemma2.1} and $Z^{c,F(Z(T_k-),\Delta\lambda
(T_k)),k}$, $k\geq0$, are independent gap processes of unregulated
evolutions without jumps with respective initial values
$F(Z(T_k-),\Delta\lambda(T_k))$, $k\geq0$. Due to the Dominated
Convergence theorem it suffices to show that the law of $Z^{z^n}(t)$
conditional on a realization of $\lambda$ converges weakly to the law
of $Z^{z^0}(t)$ conditional on the same realization of $\lambda$. But
in view of the representation (\ref{repr}), this can be shown by using
induction over the unique value of $k$ for which $t\in[T_k,T_{k+1})$
and the Feller property of $Z^{c,\cdot,k}$.

Moreover, for any $t\geq0$ we have
%
%
\begin{equation}
\lim_{n\rightarrow\infty} \frac{1}{t_n}\int_0^{t_n} (\delta_z P(s))
\,ds=\lim_{n\rightarrow\infty} \frac{1}{t_n}\int_t^{t_n+t} (\delta_z
P(s)) \,ds,
\end{equation}
where the limits are taken in the weak sense. This is a consequence of
the fact that the total variation norm of $\frac{1}{t_n}\int_0^{t_n}
(\delta_z P(s)) \,ds - \frac{1}{t_n}\int_t^{t_n+t} (\delta_z P(s))
\,ds$ is
bounded above by $\frac{2t}{t_n}$ for all $n$ with $t_n\geq t$. Hence,
by the Feller property of $Z$
\begin{eqnarray*}
\int_{W^{N-1}} f \,d(\nu P(t))&=&\lim_{n\rightarrow\infty}\int_{W^{N-1}}
P(t)f \,d\biggl(\frac{1}{t_n}\int_0^{t_n} (\delta_z P(s)) \,ds\biggr)\\
&=&\lim_{n\rightarrow\infty}\int_{W^{N-1}} f \,d\biggl(\frac
{1}{t_n}\int
_t^{t_n+t} (\delta_z P(s)) \,ds\biggr)\\
&=&\int_{W^{N-1}} f \,d\nu
\end{eqnarray*}
for all continuous bounded functions $f$. Thus, $\nu$ is an invariant
distribution of $Z$.

(4) We now prove that $\nu$ is the only invariant distribution. To this
end, consider the process $Z(n)$, $n\in\nn$. We claim that it is a
recurrent Harris chain in the sense of Section 5.6 in \cite{du}.
Indeed, in the absence of jumps the set $\{t\geq0| Z(t)\in\partial
W^{N-1}\}$ has Lebesgue measure zero almost surely and the process $Z$
evolves as a Brownian motion $(\sigma_2\beta_2(t)-\sigma_1\beta
_1(t),\ldots,\sigma_N\beta_N(t)-\sigma_1\beta_1(t))$ with constant drift
vector $(\delta_2-\delta_1,\ldots,\delta_N-\delta_1)$ in the
interior of
$W^{N-1}$ as we have seen in the proof of Lemma \ref{lemma2.1}(a). Since the
covariance matrix of the latter Brownian motion is nondegenerate and
there is a positive probability that $Z$ has no jumps in the time
interval $[0,1]$, we conclude that $Z(n)$, $n\in\nn$, is a Harris chain
on $W^{N-1}$. Noting that $\nu$ is an invariant distribution for this
chain, we conclude that the chain must be recurrent. Thus, it has a
unique invariant distribution (see \cite{du}, Section 5.6). Since
any other invariant distribution of $Z$ is an invariant distribution of
$Z(n)$, $n\in\nn$, it has to coincide with $\nu$.

(b) (1) Part (b) of the theorem can be established by using the
technique of the proof of part (a). However, we present here a softer
argument based on a monotonicity property special to the regulated
system. The latter was inspired by the proof of Lemma 1 on page 162 of
\cite{bo} which deals with a similar discrete time problem. The
crucial idea is to introduce a family of processes indexed by $\alpha
\geq0$ such that each of them evolves as the regulated system after
time $-\alpha$ and for any two indices $\alpha_1$, $\alpha_2$ the
corresponding processes are driven by the same Brownian motions and
pure jump L\'{e}vy processes after time $\max(-\alpha_1,-\alpha_2)$. To
construct such a family we start by defining auxilliary independent
standard Brownian motions $\beta_{m,1}(t),\ldots,\beta_{m,N}(t)$ and
independent i.i.d. pure jump L\'{e}vy processes $\lambda
_{m,1}(t),\ldots
,\lambda_{m,N}(t)$ of the same law as $\lambda_1(t)$ for all $m\in
\zz
_-$ on an extension of the probability space $(\Omega,\mathcal{F},\pp)$
on which $X^R$ was defined. Moreover, we let $m(t)$ be the largest
integer which is less than or equal to $t$ and $l(t)=t-m(t)$ for all
$t\in\rr$. Finally, we can use the notation to define a family of
processes $V^\alpha$ on $W^N$ indexed by $\alpha\geq0$ with the desired
properties by setting $V^\alpha(t)=b\cdot{\mathbf1}=(b,\ldots,b)$ for
$t\leq-\alpha$,
\begin{eqnarray*}
dV^\alpha_i(t)&=&\delta_i \,dt + \sigma_i \,d\beta_{m(t),i}(l(t)) +
d\lambda
_{m(t),i}(l(t))\\
&&{}+\tfrac{1}{2}\bigl(d\Lambda^\alpha_{(i-1,i)}(t)-d\Lambda^\alpha
_{(i,i+1)}(t)\bigr)\cdot1_{\{i\neq1\}}\\
&&{}+\bigl(d\Lambda^\alpha_{(0,1)}(t)-\tfrac{1}{2}\,d\Lambda^\alpha
_{(1,2)}(t)\bigr)\cdot1_{\{i=1\}}\\
&&{}+\bigl(F^R_i\bigl(V^\alpha(t-),\Delta\lambda_{m(t)}(l(t))\bigr)-\bigl(V^\alpha
_i(t-)+\Delta\lambda_{m(t),i}(l(t))\bigr)\bigr)
\end{eqnarray*}
for $-\alpha\leq t<0$ and
\begin{eqnarray*}
dV^\alpha_i(t)&=&\delta_i \,dt+\sigma_i \,d\beta_i(t)+d\lambda
_i(t)+\tfrac
{1}{2}\bigl(d\Lambda^\alpha_{(i-1,i)}(t)
-d\Lambda^\alpha_{(i,i+1)}(t)\bigr)\cdot1_{\{i\neq1\}}\\
&&{}+\bigl(d\Lambda^\alpha_{(0,1)}(t)-\tfrac{1}{2}\,d\Lambda^\alpha
_{(1,2)}(t)\bigr)\cdot1_{\{i=1\}}\\
&&{}+\bigl(F^R_i(V^\alpha(t-),\Delta\lambda(t))-\bigl(V^\alpha_i(t-)+\Delta
\lambda
_i(t)\bigr)\bigr)
\end{eqnarray*}
for $t\geq0$. Hereby, notation is as in Lemma \ref{lemma2.1}(b) with the local
times defined with respect to $V^\alpha$ instead of $Y^R$. We note that
the processes $V^\alpha$ are defined in such a way that the law of
$V^\alpha(t)$ for $t\geq-\alpha$ is the law of the ordered particle
configuration in the regulated evolution at time $t+\alpha$, started in
$b\cdot{\mathbf1}$.

(2) Moreover, we have for $V^{\alpha_1}$, $V^{\alpha_2}$ with $\alpha
_1<\alpha_2$ the inequalities
%
%
\begin{equation}
V^{\alpha_1}_1(t)\leq V^{\alpha_2}_1(t),\ldots,V^{\alpha_1}_N(t)\leq
V^{\alpha_2}_N(t)
\end{equation}
for all $t\geq-\alpha_1$. Indeed, this is clear for $t=-\alpha_1$.
Furthermore, the inequalities are preserved under the jumps of the two
processes because the jump times and the jump sizes before the
relabeling and the regulation of the particle configuration are the
same for both processes. In addition, started with each stopping time
at which a nonempty set $J\subset\{1,\ldots,N\}$ of coordinate processes
of $V^{\alpha_1}$ and $V^{\alpha_2}$ coincide, these coordinate
processes evolve in the same way until one of the particles with rank
in $J$ collides with a particle with rank not in $J$ or until there is
a jump which changes the rank of at least one particle whose rank was
originally in $J$. By distinguishing the two cases, one concludes that
the regulated evolution preserves the component-wise order $\leq$ on
configurations in $W^N$. Hence, for all $1\leq i\leq N$ we may define
%
%
\begin{equation}
V^\infty_i(t)=\uparrow\lim_{\alpha\rightarrow\infty} V^\alpha_i(t)
\end{equation}
as an element of $[b,\infty]$.

(3) Next, we claim that the random vector $V^\infty(0)$ is finite
almost surely. To this end, we observe that the process $Z^R(t)$,
$t\geq
0$, started at $Z^R(0)=0$ converges in distribution to $V^\infty
(0)-b\cdot{\mathbf1}$ for $t\rightarrow\infty$, because
%
%
\begin{equation}
Z^R(t)\stackrel{d}{=}\bigl(V^t_1(0)-b,\ldots,V^t_N(0)-b\bigr)
\end{equation}
and the right-hand side converges to $V^\infty(0)-b\cdot{\mathbf1}$ almost
surely. In addition, we have for all $1\leq i\leq N$ and all $n\in\nn
$ that
%
%
\begin{equation}
Z^R_i(n)=\sum_{j=1}^i \widetilde{Z}^R_j(n)\leq\sum_{j=1}^i D^R_j(n)
\end{equation}
with $\widetilde{Z}^R_j(t)=Z^R_j(t)-Z^R_{j-1}(t)$ for all $t\geq0$,
$1\leq j\leq N$, $Z^R_0(t)\equiv0$,
\begin{eqnarray*}
D^R_1(t)&=&\delta_1 t+\sigma_1\beta_1(t)+\sum_{k=1}^N\sum_{0\leq
s\leq
t}|\Delta\lambda_k(s)|+\Lambda^{D^R_1}(t),\\
D^R_j(t)&=&(\delta_j-\delta_{j-1})t+\sigma_j\beta_j(t)-\sigma
_{j-1}\beta
_{j-1}(t)+\sum_{k=1}^N\sum_{0\leq s\leq t}|\Delta\lambda_k(s)|
+\Lambda^{D^R_j}(t),
\end{eqnarray*}
$2\leq j\leq N$ and $\Lambda^{D^R_j}(t)$ being the local time of
$D^R_j(t)$ at zero for $1\leq j\leq N$. We now use condition (\ref
{cond2}) and follow the lines of step (2) of the proof of part (a) to
conclude that $D^R_j(t)$ converges\vspace*{1pt} in distribution to an almost surely
finite random variable in the limit $t\rightarrow\infty$ for each
$1\leq j\leq N$. Thus, for all $1\leq i\leq N$, the sequence $\sum
_{j=1}^i D^R_j(n)$, $n\in\nn$, is tight. Hence, the sequence $Z^R_i(n)$,
$n\in\nn$, is also tight, so the limit in distribution of $Z^R_i(n)$,
$n\in\nn$, which is $V^\infty_i(0)-b$, must be almost surely finite for
all $1\leq i\leq N$.

(4) Next, we prove that $Z^R$ is a Feller process. To this end, it
suffices to prove the Feller property in the absence of jumps and to
follow the lines of step (3) in the proof of part (a) to extend it to
the general case. But in the absence of jumps, the process $Z^R$ is a
reflected Brownian motion in $W^N$ in the sense of \cite{wi} as we
have shown in the proof of part (b) of the Lemma \ref{lemma2.1}. Thus, by
Theorem~1.1 and the following remark in \cite{wi}, the process $Z^R$ is a
Feller process in the absence of jumps. The Feller property and the
convergence in law of $Z^R(t)$, $t\geq0$ to $V^\infty(0)-b\cdot{\mathbf1}$
imply that the law of $V^\infty(0)-b\cdot{\mathbf1}$ is an invariant
distribution of the process $Z^R$. The uniqueness of the invariant
distribution can be deduced from the uniqueness of the invariant
distribution of the chain $Z^R(n)$, $n\in\nn$ by following step (4) of
the proof of part (a).
\end{pf*}

Now, we are able to prove Theorem \ref{theo1.3} which deals with the convergence
of the gap processes to the respective invariant distributions.
\begin{pf*}{Proof of Theorem \ref{theo1.3}}
In the course of the proof of Lemma \ref{lemma2.1}
we have seen that in the absence of jumps, the processes $Z$ and $Z^R$
are reflected Brownian motions in the sense of \cite{wi}.
Moreover, as observed in the proof of Theorem \ref{theo1.2}, their respective
covariance matrices are nondegenerate. In the presence of jumps this
shows the Harris property of the processes $Z(t)$, $t\geq0$, and
$Z^R(t)$, $t\geq0$, as well as the Harris property and the
irreducibility of the chains $Z(n)$, $n\in\nn$, and $Z^R(n)$, $n\in
\nn$.
The existence of invariant distributions (Theorem \ref{theo1.2}) implies that the
processes $Z(t)$, $t\geq0$, and $Z^R(t)$, $t\geq0$, are positive Harris
recurrent in the sense of \cite{mt}. Thus, Theorem 6.1 of \cite{mt}
shows that $Z(t)$, $t\geq0$, and $Z^R(t),t\geq0$, converge in total
variation to their respective invariant distributions.
\end{pf*}

\section{Systems of infinitely many particles}\label{sec3}

\subsection{Existence and uniqueness in law of the processes}\label{sec31}

From now on we let $I=\nn$. Due to the fact that the processes
$L_1(t),L_2(t),\ldots$ are i.i.d., an application of the second
Borel--Cantelli lemma shows immediately that in both evolutions
infinitely many jumps of the particles will occur almost surely on each
nonempty time interval, if $L_1(t)$ is not identically equal to zero.
Hence, the proof of existence of the finite particle systems does not
carry over to the infinite case. Also, the proof in \cite{pp} for
the case of the unregulated system without jumps, which uses Girsanov's
theorem for Brownian motion, cannot be applied here, since it would
prove the existence of the solution to the unregulated version of (\ref
{sde1}) in the case of a certain nontrivial dependence structure
between the processes $B_1(t),B_2(t),\ldots$ and
$L_1(t),L_2(t),\ldots.$
Instead, we prove the existence and the uniqueness in law of the
infinite particle systems by a bound on the tail of the distribution of
the running supremum of an integrable L\'{e}vy process.
\begin{prop}\label{prop3.1}
If the initial configuration $X_1(0),X_2(0),\ldots$ of the
particles satisfies
%
%
\begin{equation}\label{supcond}
\sum_{i=1}^\infty\pp\Bigl(\sup_{0\leq s\leq t} \bigl(-\sigma_1
B_1(s)-L_1(s)\bigr)\geq X_i(0)-y\Bigr) < \infty
\end{equation}
for all $t\geq0$ and $y\in\rr$, then the unregulated version of
(\ref{sde1}) has a unique weak solution and the corresponding ordered
particle system $Y_1(t),Y_2(t),\ldots$ is well defined for all $t\geq
0$.

In particular, condition (\ref{supcond}) is satisfied if there are
constants $\gamma_1>0$, $\gamma_2\in\rr$ with
%
%
\begin{equation}\label{lincond}
X_i(0)\geq\gamma_1 i+\gamma_2, \qquad i\geq1.
\end{equation}
The same statements are true for the regulated system (\ref{sde1}),
(\ref{sde2}).
\end{prop}
\begin{pf}
(1) We prove the proposition only for the unregulated system,
since the assertion for the regulated system can be shown in the same
way by setting $M=1$. To this end, we assume (\ref{supcond}) and
introduce the probability space
\[
\biggl(\prod_{(N,x)} \Omega^{(N,x)}\times\Omega^{B,L},\bigotimes_{(N,x)}
\mathcal{F}^{(N,x)}\otimes\mathcal{F}^{B,L},\bigotimes_{(N,x)} \pp
^{(N,x)}\otimes\pp^{B,L}\biggr).
\]
Hereby, the products are over $\{(N,x)|N\geq1,x\in\rr^N,x_1\leq
\cdots
\leq x_N\}$, on each $(\Omega^{(N,x)},\mathcal{F}^{(N,x)},\pp^{(N,x)})$
the unregulated system $X^{N,x_1,\ldots,x_N}$ with $N$ particles,
initial configuration $x_1\leq x_2\leq\cdots\leq x_N$ and parameters
$\delta_1,\ldots,\delta_N,\sigma_1,\ldots,\sigma_1$ are defined and
$(\Omega^{B,L},\mathcal{F}^{B,L},\pp^{B,L})$ is a probability space on
which the i.i.d. Brownian motions $B_1(t),B_2(t),\ldots$ and the i.i.d.
pure jump L\'{e}vy processes $L_1(t),L_2(t),\ldots$ independent of the
Brownian motions are defined. We call the product space above $(\Omega
^\infty,\mathcal{F}^\infty,\pp^\infty)$. Next, we recall the definition
of the constant $M$ introduced in Assumption \ref{assump1.1}, set $G_i(t)=\sigma_1
B_i(t)+L_i(t)$ for all $t\geq0$ and $i\in I$ for the sake of shorter
notation and define the process $X$, the sets $A_0\subset A_1\subset
\cdots$ and the stopping times $0=\rho_0\leq\rho_1\leq\cdots$
inductively by
\begin{eqnarray*}
A_k&=&\bigl\{i\geq1|\exists1\leq j\leq M, 0\leq s\leq\rho_k\dvtx
X_i(s)=X_{(j)}(s)\bigr\},\\
X_i(s)&=&X_i^{|A_k|,X_{(1)}(\rho_k),\ldots,X_{(|A_k|)}(\rho
_k)}(s-\rho
_k),\qquad i\in A_k, \rho_k\leq s<\rho_{k+1},\\
X_i(s)&=&X_i(\rho_k)+\delta_M(s-\rho_k)+G_i(s)-G_i(\rho_k),\qquad
i\notin
A_k, \rho_k\leq s<\rho_{k+1},\\
\rho_{k+1}&=&\inf\bigl\{s\geq\rho_k|\exists1\leq j\leq M, i\notin A_k\dvtx
X_i(s)=X_{(j)}(s)\bigr\}.
\end{eqnarray*}
Fixing a $t>0$ we observe that by condition (\ref{supcond}) and the
first Borel--Cantelli lemma there are almost surely only finitely many
particles which visit the interval $(-\infty,\widetilde{y}]$ up to time
$t$ for any fixed $\widetilde{y}\in\rr$. This follows by setting
$y=\widetilde{y}+{\max_{i\geq1}} |\delta_i|\cdot t$ in (\ref{supcond})
and recalling that the set of collision times of distinct particles in
the finite unregulated particle system has Lebesgue measure zero almost
surely (see the proof of Lemma \ref{lemma2.1}). Furthermore, we note that by
choosing a large enough $\widetilde{y}$ we can make
\[
\sum_{j=1}^M \pp\Bigl(X_j(0)+\sup_{0\leq s\leq t} \bigl(\sigma_1 B_1(s) +
L_1(s)\bigr) + \max_{i\geq1} |\delta_i|\cdot t>\widetilde{y}\Bigr)
\]
as small as we want. We also note that in the regulated case this
expression should be replaced by
\[
\pp\Bigl(\sup_{0\leq s\leq t} \bigl(X_1(0)+\delta_1 s + \sigma_1 B_1(s) +
L_1(s) + \Lambda_{(1,b)}(s)\bigr)>\widetilde{y}\Bigr).
\]
This observation implies that there exists a $K=K(\omega)$ such that
$K<\infty$, $\rho_K\leq t<\rho_{K+1}$ and the sets $A_0,\ldots,A_K$ are
finite almost surely. Thus, $X$ is well defined on $[0,t]$ for almost
every $\omega\in\Omega^\infty$. Since $t>0$ was arbitrary, we have
shown that $X$ is well defined on $[0,\infty)$ for almost every
$\omega
\in\Omega^\infty$. By its construction, $X$ is a weak solution to the
unregulated version of the system (\ref{sde1}). Moreover, at any time
$t\geq0$ the ordered particle system $Y_1(t),Y_2(t),\ldots$ is
well defined, since there are finitely many particles on each interval
of the form $(-\infty,\widetilde{y}]$ almost surely [due to (\ref
{supcond}) with $y=\widetilde{y}+{\max_{i\geq1}} |\delta_i|\cdot t$ and
the first Borel--Cantelli lemma].

(2) For the uniqueness part let $X'$ be a different weak solution to
the unregulated version of (\ref{sde1}). Then we can define inductively
the sets $A'_0\subset A'_1\subset\cdots$ and the stopping times
$0=\rho
'_0\leq\rho'_1\leq\cdots$ by
%
%
\begin{eqnarray}
A'_k&=&\bigl\{i\geq1|\exists1\leq j\leq M, 0\leq s\leq\rho'_k\dvtx
X'_i(s)=X'_{(j)}(s)\bigr\},\\
\rho'_{k+1}&=&\inf\bigl\{s\geq\rho'_k|\exists1\leq j\leq M, i\notin A_k\dvtx
X'_i(s)=X'_{(j)}(s)\bigr\}.
\end{eqnarray}
Due to the uniqueness of the weak solution to the unregulated system
(\ref{sde1}) in the case of finitely many particles, the joint
distribution of $\rho'_0,\rho'_1,\ldots$ and $X'$ on $[\rho'_0,\rho
'_1],[\rho'_1,\rho'_2],\ldots$ has to coincide with the joint
distribution of $\rho_0,\rho_1,\ldots$ and $X$ on $[\rho_0,\rho
_1],[\rho
_1,\rho_2],\ldots.$ Thus, the law of $X'$ is the same as the law of
$X$.

(3) To prove that (\ref{lincond}) implies (\ref{supcond}) we first
observe that
%
%
\begin{equation}
\ev\bigl[\max\bigl(|\sigma_1 B_1(t)+L_1(t)|,1\bigr)\bigr]<\infty
\end{equation}
for any fixed $t\geq0$. Noting that $\max(|x|,1)$ is a nonnegative
continuous submultiplicative function (see Proposition 25.4 of \cite
{sa}), we conclude from the Theorem~25.18 of \cite{sa} that
%
%
\begin{equation}
\ev\Bigl[\sup_{0\leq s\leq t}|\sigma_1 B_1(s)+L_1(s)|\Bigr]<\infty
\end{equation}
for any fixed $t\geq0$. Thus, (\ref{lincond}) implies (\ref
{supcond}).
\end{pf}

\subsection{Tightness of the infinite regulated system}\label{sec32}

We can now prove Proposition~\ref{prop1.4} which guarantees the tightness of the
gap process in the infinite regulated evolution.
\begin{pf*}{Proof of Proposition \ref{prop1.4}}
We start with the observation that
for each $i\in I$ the process $X_i^R(t)$, $t\geq0$ and the
corresponding process of regulations $R_i(t)$, $t\geq0$, solve the
system of stochastic differential equations
\begin{eqnarray*}
dX^R_i(t)&=&\delta_1 \,dt + \sigma_1 \,dB_i(t) + dL_i(t)+ dR_i(t),\\
dR_i(t)&=&\bigl(b-X^R_i(t-)-\Delta L_i(t)\bigr)1_{\{X^R_i(t-)+\Delta L_i(t)<b\}
}+d\Lambda_{(i,b)}(t).
\end{eqnarray*}
This is due to the fact that the drift and the volatility sequences are
constant and that the set of collision times of distinct particles has
Lebesgue measure zero almost surely. The latter statement is a
consequence of the corresponding property of the regulated system with
finitely many particles [see the proof of part (b) of Lemma \ref{lemma2.1}] and
the construction of the infinite particle systems (see the proof of
Proposition \ref{prop3.1}). Next, for each $i\in I$ we introduce the process
%
%
\begin{equation}
H^R_i(t)=X^R_i(0)+\delta_1 t + \sigma_1 B_i(t) + \sum_{0\leq s\leq
t}|\Delta L_i(s)| + \Lambda^{H^R_i}(t),
\end{equation}
where $\Lambda^{H^R_i}(t)$ is the local time of $H^R_i$ at $b$. Using
condition (\ref{cond3}) and arguing as in the steps\vspace*{1pt} (1) and (2) of the
proof of Theorem \ref{theo1.2}(a), we conclude that $X^R_i(t)\leq H^R_i(t)$ for
all $t\geq0$ and $i\in I$ almost surely and, in addition, that for each
$i\in I$ the process $H^R_i(t)$, $t\geq0$, converges in law to an almost
surely finite random variable $H^R_i(\infty)$ whose law does not depend
on $X^R_i(0)$ and $i$. Moreover, due to the independence of the
processes $H^R_i$, $i\in I$, the random vector $(H^R_1(t),\ldots
,H^R_j(t))$ converges in distribution to $(H^R_1(\infty),\ldots
,H^R_j(\infty))$ for any $j\in I$ where $H^R_i(\infty)$, $i\geq2$, are
chosen as independent copies of $H^R_1(\infty)$. Finally, the chain of
inequalities
%
%
\begin{equation}
0\leq Z^R_j(t)\leq\max_{1\leq i\leq j} X^R_i(t)-b\leq\max_{1\leq
i\leq
j} H^R_i(t)-b
\end{equation}
for all $t\geq0$ shows that the family $Z^R_j(t)$, $t\geq0$, is tight
for all $j\in I$. This yields the tightness of $Z^R(t)$, $t\geq0$, on
$\rr_+^\nn$ with the product topology.
\end{pf*}

\section*{Acknowledgments}
The author thanks Amir Dembo for his invaluable comments and
suggestions throughout the preparation of this work. He is also
grateful to George Papanicolaou and Jim Pitman for helpful discussions
and to two anonymous referees for a careful reading of the paper and
their comments.

%

%
\printaddresses


\begin{thebibliography}{34}

\bibitem{aa}
%
\begin{barticle}[mr]
\bauthor{\bsnm{Arguin},~\bfnm{Louis-Pierre}\binits{L.-P.}} \AND
\bauthor{\bsnm{Aizenman},~\bfnm{Michael}\binits{M.}}
(\byear{2009}).
\btitle{On the structure of quasi-stationary competing particle systems}.
\bjournal{Ann. Probab.}
\bvolume{37}
\bpages{1080--1113}.
\bid{doi={10.1214/08-AOP429}, mr={2537550}}
\end{barticle}
%
\endbibitem

\bibitem{as}
%
\begin{bbook}[mr]
\bauthor{\bsnm{Asmussen},~\bfnm{S{\o}ren}\binits{S.}}
(\byear{2003}).
\btitle{Applied Probability and Queues: Stochastic Modelling and Applied Probability},
\bedition{2nd} ed.
\bseries{Applications of Mathematics (New York)}
\bvolume{51}.
\bpublisher{Springer}, \baddress{New York}.
\bid{mr={1978607}}
\end{bbook}
%
\endbibitem

\bibitem{bf}
%
\begin{barticle}[mr]
\bauthor{\bsnm{Banner},~\bfnm{Adrian~D.}\binits{A.~D.}},
\bauthor{\bsnm{Fernholz},~\bfnm{Robert}\binits{R.}} \AND
\bauthor{\bsnm{Karatzas},~\bfnm{Ioannis}\binits{I.}}
(\byear{2005}).
\btitle{Atlas models of equity markets}.
\bjournal{Ann. Appl. Probab.}
\bvolume{15}
\bpages{2296--2330}.
\bid{doi={10.1214/105051605000000449}, mr={2187296}}
\end{barticle}
%
\endbibitem

\bibitem{bg}
%
\begin{barticle}[mr]
\bauthor{\bsnm{Banner},~\bfnm{Adrian~D.}\binits{A.~D.}} \AND
\bauthor{\bsnm{Ghomrasni},~\bfnm{Raouf}\binits{R.}}
(\byear{2008}).
\btitle{Local times of ranked continuous semimartingales}.
\bjournal{Stochastic Process. Appl.}
\bvolume{118}
\bpages{1244--1253}.
\bid{doi={10.1016/j.spa.2007.08.001}, mr={2428716}}
\end{barticle}
%
\endbibitem

\bibitem{ba}
%
\begin{barticle}[mr]
\bauthor{\bsnm{Bardhan},~\bfnm{Indrajit}\binits{I.}}
(\byear{1995}).
\btitle{Further applications of a general rate conservation law}.
\bjournal{Stochastic Process. Appl.}
\bvolume{60}
\bpages{113--130}.
\bid{doi={10.1016/0304-4149(95)00052-6}, mr={1362322}}
\end{barticle}
%
\endbibitem\vadjust{\goodbreak}

\bibitem{bp}
%
\begin{barticle}[mr]
\bauthor{\bsnm{Bass},~\bfnm{R.~F.}\binits{R.~F.}} \AND
\bauthor{\bsnm{Pardoux},~\bfnm{{\'E}.}\binits{{\'E}.}}
(\byear{1987}).
\btitle{Uniqueness for diffusions with piecewise constant coefficients}.
\bjournal{Probab. Theory Related Fields}
\bvolume{76}
\bpages{557--572}.
\bid{doi={10.1007/BF00960074}, mr={0917679}}
\end{barticle}
%
\endbibitem



\bibitem{bo}
%
\begin{bbook}[mr]
\bauthor{\bsnm{Borovkov},~\bfnm{A.~A.}\binits{A.~A.}}
(\byear{1976}).
\btitle{Stochastic Processes in Queueing Theory}.
\bpublisher{Springer}, \baddress{New York}.
\bid{mr={0391297}}
\end{bbook}
%
\endbibitem

\bibitem{cp}
%
\begin{bmisc}[auto:STB|2010-11-18|09:18:59]
\bauthor{\bsnm{Chatterjee},~\bfnm{S.}\binits{S.}} \AND
\bauthor{\bsnm{Pal},~\bfnm{S.}\binits{S.}}
(\byear{2009}).
\bhowpublished{A phase transition behaviour for Brownian motions interacting through
their ranks. Available at}
\href{http://arxiv.org/abs/arXiv:0706.3558v2}{arXiv:0706.3558v2}.
\end{bmisc}
%
\endbibitem


\bibitem{ct}
%
\begin{bbook}[mr]
\bauthor{\bsnm{Cont},~\bfnm{Rama}\binits{R.}} \AND
\bauthor{\bsnm{Tankov},~\bfnm{Peter}\binits{P.}}
(\byear{2004}).
\btitle{Financial Modelling with Jump Processes}.
\bpublisher{Chapman and Hall/CRC}, \baddress{Boca Raton, FL}.
\bid{mr={2042661}}
\end{bbook}
%
\endbibitem


\bibitem{du}
%
\begin{bbook}[mr]
\bauthor{\bsnm{Durrett},~\bfnm{Richard}\binits{R.}}
(\byear{1996}).
\btitle{Probability: Theory and Examples}, \bedition{2nd} ed.
\bpublisher{Duxbury Press}, \baddress{Belmont, CA}.
\bid{mr={1609153}}
\end{bbook}
%
\endbibitem

\bibitem{fe}
%
\begin{bbook}[mr]
\bauthor{\bsnm{Fernholz},~\bfnm{E.~Robert}\binits{E.~R.}}
(\byear{2002}).
\btitle{Stochastic Portfolio Theory: Stochastic Modelling and Applied Probability}.
\bseries{Applications of Mathematics (New York)}
\bvolume{48}.
\bpublisher{Springer}, \baddress{New York}.
\bid{mr={1894767}}
\end{bbook}
%
\endbibitem

\bibitem{ha}
%
\begin{bincollection}[mr]
\bauthor{\bsnm{Harrison},~\bfnm{J.~Michael}\binits{J.~M.}}
(\byear{1988}).
\btitle{Brownian models of queueing networks with heterogeneous customer
populations}.
In \bbooktitle{Stochastic Differential Systems, Stochastic Control
Theory and
Applications ({M}inneapolis, {M}inn., 1986)}.
\bseries{The IMA Volumes in Mathematics and its Applications}
\bvolume{10}
\bpages{147--186}.
\bpublisher{Springer}, \baddress{New York}.
\bid{mr={0934722}}
\end{bincollection}
%
\endbibitem

\bibitem{hn}
%
\begin{barticle}[mr]
\bauthor{\bsnm{Harrison},~\bfnm{J.~Michael}\binits{J.~M.}} \AND
\bauthor{\bsnm{Nguyen},~\bfnm{Vi{\^e}n}\binits{V.}}
(\byear{1993}).
\btitle{Brownian models of multiclass queueing networks: Current
status and
open problems}.
\bjournal{Queueing Syst.}
\bvolume{13}
\bpages{5--40}.
\bid{mr={1218842}}
\end{barticle}
%
\endbibitem

\bibitem{ip}
%
\begin{barticle}[auto:STB|2010-11-18|09:18:59]
\bauthor{\bsnm{Ichiba},~\bfnm{T.}\binits{T.}},
\bauthor{\bsnm{Papathanakos},~\bfnm{V.}\binits{V.}},
\bauthor{\bsnm{Banner},~\bfnm{A.}\binits{A.}},
\bauthor{\bsnm{Karatzas},~\bfnm{I.}\binits{I.}} \AND
\bauthor{\bsnm{Fernholz},~\bfnm{R.}\binits{R.}}
(\byear{2011}).
\btitle{Hybrid atlas models}.
\bjournal{Ann. Appl. Probab.}
\bvolume{21}
\bpages{609--644}.
\end{barticle}
%
\endbibitem

\bibitem{ka}
%
\begin{bbook}[mr]
\bauthor{\bsnm{Kallenberg},~\bfnm{Olav}\binits{O.}}
(\byear{2002}).
\btitle{Foundations of Modern Probability},
\bedition{2nd} ed.
\bpublisher{Springer}, \baddress{New York}.
\bid{mr={1876169}}
\end{bbook}
%
\endbibitem


\bibitem{ms}
%
\begin{barticle}[mr]
\bauthor{\bsnm{McKean},~\bfnm{H.~P.}\binits{H.~P.}} \AND
\bauthor{\bsnm{Shepp},~\bfnm{L.~A.}\binits{L.~A.}}
(\byear{2005}).
\btitle{The advantage of capitalism vs. socialism depends on the criterion}.
\bjournal{Zap. Nauchn. Sem. S.-Peterburg. Otdel. Mat. Inst. Steklov. (POMI)}
\bvolume{328}
\bpages{160--168, 279--280}.
\bid{doi={10.1007/s10958-006-0374-5}, mr={2214539}}
\end{barticle}
%
\endbibitem

\bibitem{me}
%
\begin{barticle}[auto:STB|2010-11-18|09:18:59]
\bauthor{\bsnm{Merton},~\bfnm{R.~C.}\binits{R.~C.}}
(\byear{1976}).
\btitle{Option pricing when underlying stock returns are discontinuous}.
\bjournal{Journal of Financial Economics}
\bvolume{3}
\bpages{125--144}.
\end{barticle}
%
\endbibitem

\bibitem{mt}
%
\begin{barticle}[mr]
\bauthor{\bsnm{Meyn},~\bfnm{Sean~P.}\binits{S.~P.}} \AND
\bauthor{\bsnm{Tweedie},~\bfnm{R.~L.}\binits{R.~L.}}
(\byear{1993}).
\btitle{Stability of {M}arkovian processes. {II}. {C}ontinuous-time processes
and sampled chains}.
\bjournal{Adv. in Appl. Probab.}
\bvolume{25}
\bpages{487--517}.
\bid{doi={10.2307/1427521}, mr={1234294}}
\end{barticle}
%
\endbibitem

\bibitem{pp}
%
\begin{barticle}[mr]
\bauthor{\bsnm{Pal},~\bfnm{Soumik}\binits{S.}} \AND
\bauthor{\bsnm{Pitman},~\bfnm{Jim}\binits{J.}}
(\byear{2008}).
\btitle{One-dimensional {B}rownian particle systems with rank-dependent
drifts}.
\bjournal{Ann. Appl. Probab.}
\bvolume{18}
\bpages{2179--2207}.
\bid{doi={10.1214/08-AAP516}, mr={2473654}}
\end{barticle}
%
\endbibitem

\bibitem{rw}
%
\begin{barticle}[mr]
\bauthor{\bsnm{Reiman},~\bfnm{M.~I.}\binits{M.~I.}} \AND
\bauthor{\bsnm{Williams},~\bfnm{R.~J.}\binits{R.~J.}}
(\byear{1988}).
\btitle{A boundary property of semimartingale reflecting {B}rownian motions}.
\bjournal{Probab. Theory Related Fields}
\bvolume{77}
\bpages{87--97}.
\bid{doi={10.1007/BF01848132}, mr={0921820}}
\end{barticle}
%
\endbibitem

\bibitem{ry}
%
\begin{bbook}[mr]
\bauthor{\bsnm{Revuz},~\bfnm{Daniel}\binits{D.}} \AND
\bauthor{\bsnm{Yor},~\bfnm{Marc}\binits{M.}}
(\byear{1999}).
\btitle{Continuous Martingales and {B}rownian Motion},
\bedition{3rd} ed.
\bseries{Grundlehren der Mathematischen Wissenschaften [Fundamental Principles
of Mathematical Sciences]}
\bvolume{293}.
\bpublisher{Springer}, \baddress{Berlin}.
\bid{mr={1725357}}
\end{bbook}
%
\endbibitem

\bibitem{ra}
%
\begin{barticle}[mr]
\bauthor{\bsnm{Ruzmaikina},~\bfnm{Anastasia}\binits{A.}} \AND
\bauthor{\bsnm{Aizenman},~\bfnm{Michael}\binits{M.}}
(\byear{2005}).
\btitle{Characterization of invariant measures at the leading edge for
competing particle systems}.
\bjournal{Ann. Probab.}
\bvolume{33}
\bpages{82--113}.
\bid{doi={10.1214/009117904000000865}, mr={2118860}}
\end{barticle}
%
\endbibitem

\bibitem{sa}
%
\begin{bbook}[mr]
\bauthor{\bsnm{Sato},~\bfnm{Ken-Iti}\binits{K.-I.}}
(\byear{1999}).
\btitle{L\'evy Processes and Infinitely Divisible Distributions}.
\bseries{Cambridge Studies in Advanced Mathematics}
\bvolume{68}.
\bpublisher{Cambridge Univ. Press}, \baddress{Cambridge}.
\bid{mr={1739520}}
\end{bbook}
%
\endbibitem

\bibitem{sh}
%
\begin{barticle}[mr]
\bauthor{\bsnm{Shkolnikov},~\bfnm{Mykhaylo}\binits{M.}}
(\byear{2009}).
\btitle{Competing particle systems evolving by i.i.d. increments}.
\bjournal{Electron. J. Probab.}
\bvolume{14}
\bpages{728--751}.
\bid{mr={2486819}}
\end{barticle}
%
\endbibitem

%

\bibitem{tw}
%
\begin{barticle}[mr]
\bauthor{\bsnm{Taylor},~\bfnm{L.~M.}\binits{L.~M.}} \AND
\bauthor{\bsnm{Williams},~\bfnm{R.~J.}\binits{R.~J.}}
(\byear{1993}).
\btitle{Existence and uniqueness of semimartingale reflecting {B}rownian
motions in an orthant}.
\bjournal{Probab. Theory Related Fields}
\bvolume{96}
\bpages{283--317}.
\bid{doi={10.1007/BF01292674}, mr={1231926}}
\end{barticle}
%
\endbibitem


\bibitem{wi}
%
\begin{barticle}[mr]
\bauthor{\bsnm{Williams},~\bfnm{R.~J.}\binits{R.~J.}}
(\byear{1987}).
\btitle{Reflected {B}rownian motion with skew symmetric data in a polyhedral
domain}.
\bjournal{Probab. Theory Related Fields}
\bvolume{75}
\bpages{459--485}.
\bid{doi={10.1007/BF00320328}, mr={0894900}}
\end{barticle}
%
\endbibitem\vadjust{\goodbreak}

\bibitem{wi2}
%
\begin{bincollection}[mr]
\bauthor{\bsnm{Williams},~\bfnm{R.~J.}\binits{R.~J.}}
(\byear{1995}).
\btitle{Semimartingale reflecting {B}rownian motions in the orthant}.
In \bbooktitle{Stochastic Networks}.
\bseries{The IMA Volumes in Mathematics and its Applications}
\bvolume{71}
\bpages{125--137}.
\bpublisher{Springer}, \baddress{New York}.
\bid{mr={1381009}}
\end{bincollection}
%
\endbibitem



\end{thebibliography}
\end{document}